\documentclass[12pt]{article}
\usepackage{amssymb,amsmath}
\newcommand{\f}{\mathbf{f}}
\newcommand{\g}{\mathbf{g}}

\renewcommand{\u}{\mathbf{u}}
\newcommand{\Z}{{\mathbb Z}}

\newcommand{\Fp}{{\mathbb F}_p}
\newcommand{\Q}{{\mathbb Q}}
\newcommand{\R}{\mathcal{R}}
\renewcommand{\v}{\vec{v}}
\newcommand{\N}{{\mathcal N}}
\newcommand{\kmod}[1]{\:\:(\text{mod }#1)}
\newcommand{\M}{{\rm M}}
\newcommand{\cbar}{\overline{c}}
\newcommand{\qbar}{\overline{q}}
\newcommand{\Qbar}{\overline{Q}}
\newcommand{\Pbar}{\overline{P}}
\newcommand{\Nbar}{\overline{\mathcal N}}
\newcommand{\fbar}{\overline{\mathbf{f}}}
\newcommand{\gbar}{\overline{\mathbf{g}}}
\newcommand{\phibar}{\overline{\phi}}
\newcommand{\ubar}{\overline{\mathbf{u}}}
\renewcommand{\hbar}{\overline{\mathbf{h}}}
\newcommand{\proof}{\noindent{\em Proof: }}
\newcommand{\qed}{\hspace{\fill}$\square$}
\newcommand{\ra}{\rightarrow}

\newcommand{\dst}{\displaystyle}

\newtheorem{theorem}{Theorem}
\newtheorem{lemma}[theorem]{Lemma}
\newtheorem{prop}[theorem]{Proposition}
\newtheorem{cor}[theorem]{Corollary}
\newtheorem{case}{Case}
\newenvironment{remark}{\noindent\refstepcounter{theorem}{\bf Remark \arabic{theorem}} }{}

\title{How close are $p$th powers in the Nottingham group?}
\author{Kevin Keating \\
Department of Mathematics \\
University of Florida \\
Gainesville, FL 32611 \\
USA \\[.2cm]
{\tt keating@math.ufl.edu}}
\date{}
\begin{document}

\maketitle

\begin{abstract}
\noindent
Let $F$ be a field of characteristic $p>0$ and let $f,g$ be
elements of the Nottingham group $\N(F)$ such that $f$ has
depth $k$ and $gf^{-1}$ has depth $n\ge k$.  We find the best
possible lower bound for the depth of $g^pf^{-p}$.
\end{abstract}

     Let $R$ be a commutative ring with 1 which has
characteristic $p>0$, and let $\N=\N(R)$
denote the Nottingham group over $R$.  Thus $\N$ is the set of
all formal power series $f(x)\in R[[x]]$ with leading term $x$,
and the product of $f,g\in\N$ is defined to be
$(fg)(x)=f(g(x))$.  For each $k\geq1$ define a normal subgroup
$\N_k\trianglelefteq\N$ by setting
\begin{equation}
\N_k=\{f\in\N:f(x)\equiv x\kmod{x^{k+1}}\}.
\end{equation}
The depth of $f\in\N$ is defined to be $D(f)=\sup\{k:f\in\N_k\}$.

     Let $n\ge k\ge1$ and let $k_0$ be the least nonnegative
residue of $k$ modulo $p$.  We define a nonnegative integer
$e(k,n)$ as follows:
\begin{equation} \label{edef}
e(k,n)=
\begin{cases}
0&\text{ if $p\mid k$ and $n=k$,} \\
1&\text{ if $p\mid k$, $p\mid n$, and $n>k$,} \\
0&\text{ if $p\mid k$ and $p\nmid n$,} \\
i&\text{ if $p\nmid k$ and $n\equiv2k-i\kmod{p}$ for some
$0\le i\le k_0$,} \\
k_0&\text{ if $p\nmid k$ and $n\not\equiv2k-i\kmod{p}$ for all
$0\le i\le k_0$.}
\end{cases}
\end{equation}
Our main result is the following:

\begin{theorem} \label{main}
(a) Let $f$ and $g$ be elements of $\N(R)$ such that $D(f)\ge k$
and $D(gf^{-1})\ge n$.  Then
\begin{equation} \label{m}
D(g^pf^{-p})\ge n+(p-1)k+e(k,n).
\end{equation}
(b) There exist $f,g\in\N(R)$ such that $D(f)=k$,
$D(gf^{-1})=n$, and 
\begin{equation}
D(g^pf^{-p})=n+(p-1)k+e(k,n).
\end{equation}
\end{theorem}

     The following corollary generalizes Theorem~\ref{main}(a)
to higher powers of $p$.  It would be interesting to know
whether the bound given here is the best possible.

\begin{cor} \label{pmcor}
Let $f,g\in\N(R)$ be such that $D(f)\ge k$ and $D(gf^{-1})\ge n$.
Then for all $m\ge1$ we have
\begin{equation} \label{pm}
D(g^{p^m}f^{-p^m})\ge n+(p^m-1)k+\frac{p^m-p}{p-1}k_0+e(k,n).
\end{equation}
\end{cor}

\proof By repeated application of Lemma~\ref{powers} below
we get
\begin{equation} \label{Dgpm}
D(f^{p^i})\ge p^ik+\frac{p^i-1}{p-1}k_0
\end{equation}
for all $i\ge1$.  It follows from (\ref{edef}) that
\begin{equation} \label{nnp}
n+(p-1)k+e(k,n)\equiv
\begin{cases}
n\pmod{p}&\text{if $e(k,n)=k_0$,} \\
k\pmod{p}&\text{if $0\le e(k,n)<k_0$,} \\
1\pmod{p}&\text{if $p\mid k$, $p\mid n$, and $n>k$.}
\end{cases}
\end{equation}
Using (\ref{m}), (\ref{Dgpm}), and (\ref{nnp})
we can iteratively compute lower bounds for $D(g^{p^i}f^{-p^i})$.
For $i\geq1$ we get $D(g^{p^i}f^{-p^i})\geq d_i$, where
$d_1=n+(p-1)k+e(k,n)$ and
\begin{equation}
d_{i+1}=d_i+(p-1)\left(p^ik+\frac{p^i-1}{p-1}k_0\right)+k_0.
\end{equation}
By summing the terms we get
\begin{equation}
d_m=n+(p^m-1)k+\frac{p^m-p}{p-1}k_0+e(k,n),
\end{equation}
as required. \qed \medskip

     Let $S$ be a commutative ring with 1 which has
characteristic $p$ and let
$\sigma:R\ra S$ be a unitary ring homomorphism.  Then $\sigma$
induces a group homomorphism $\N(R)\ra\N(S)$ which we denote by
$f\mapsto f^{\sigma}$.  We clearly have
$D(f^{\sigma})\ge D(f)$.  Let
$\R=\Fp[r_k,r_{k+1},\dots, s_n,s_{n+1},\dots]$,
where $\Fp=\Z/p\Z$ is the field with $p$ elements and
$r_i,s_j$ are variables.  Also set
\begin{align}
\f(x)&=x+r_kx^{k+1}+r_{k+1}x^{k+2}+\cdots \\
\u_1(x)&=x+s_nx^{n+1}+s_{n+1}x^{n+2}+\cdots
\end{align}
and $\g=\u_1\f$.  Then $\f$ is a generic element of $\N(\R)$
of depth $k$, and $\u_1$ is a generic element of $\N(\R)$ of
depth $n$.

     Let $f,g\in\N(R)$ satisfy $D(f)\ge k$ and
$D(gf^{-1})\ge n$.  Then we have
\begin{align}
f(x)&=x+a_kx^{k+1}+a_{k+1}x^{k+2}+\cdots \\
gf^{-1}(x)&=x+b_nx^{n+1}+b_{n+1}x^{n+2}+\cdots
\end{align}
with $a_i,b_j\in R$.
Let $\sigma:\R\ra R$ be the unique homomorphism such that
$\sigma(r_i)=a_i$ for $i\ge k$ and $\sigma(s_j)=b_j$ for
$j\ge n$.  Then $\f^{\sigma}=f$,
$(\g\f^{-1})^{\sigma}=gf^{-1}$, and hence $\g^{\sigma}=g$.
Therefore to prove Theorem~\ref{main}(a) it suffices to show
\begin{equation}
D(\g^p\f^{-p})\ge n+(p-1)k+e(k,n).
\end{equation}
Since $\Fp$ is a subring of every ring of characteristic
$p$, it suffices to prove Theorem~\ref{main}(b) in the case
$R=\Fp$.  We define a specialization to be a homomorphism
$\sigma:\R\ra\Fp$.  Associated to a specialization $\sigma$
we have elements $\f^{\sigma}$, $\g^{\sigma}$, $\u_1^{\sigma}$
of $\N(\Fp)$. \medskip

\begin{remark} \label{interp}
\rm Theorem~\ref{main} can be expressed entirely in terms of the
generic power series $\f(x)$.  Theorem~\ref{main}(a) is
equivalent to the statement that for all
$i\le n+(p-1)k+e(k,n)$ the coefficient of $x^i$ in $\f^p(x)$
does not depend on any $r_j$ with $j\ge n$.
Theorem~\ref{main}(b) is equivalent to the statement that
there exist specializations $\sigma,\tau$ such that
$\sigma(r_j)=\tau(r_j)$ for $k\le j<n$,
$\sigma(r_n)\not=\tau(r_n)$, and
$D((\f^{\tau})^p(\f^{\sigma})^{-p})=n+(p-1)k+e(k,n)$. \qed
\end{remark} \medskip

     The following lemma is useful in the proof of
Theorem~\ref{main}(b):

\begin{lemma} \label{shift}
Suppose there are $n'>n\ge k$ such that
Theorem~\ref{main}(a) holds for $(k,n)$,
Theorem~\ref{main}(b) holds for $(k,n')$, and
\begin{equation}
n+(p-1)k+e(k,n)=n'+(p-1)k+e(k,n').
\end{equation}
Then Theorem~\ref{main}(b) holds for $(k,n)$.
\end{lemma}

\proof Since Theorem~\ref{main}(b) holds for $(k,n')$ there
are $f,g\in\N(R)$ such that $D(f)=k$, $D(gf^{-1})=n'$, and
$D(g^pf^{-p})=n+(p-1)k+e(k,n)$.  Choose $h\in\N(R)$ such that
$D(hf^{-1})=n$.  Then we have
\begin{equation}
D(h^pf^{-p})\ge D(\g^p\f^{-p})\ge n+(p-1)k+e(k,n).
\end{equation}
If $D(h^pf^{-p})=n+(p-1)k+e(k,n)$ then $f$, $h$ satisfy the
conditions of
Theorem~\ref{main}(b).  If $D(h^pf^{-p})>n+(p-1)k+e(k,n)$ then
${D(h^pg^{-p})=n+(p-1)k+e(k,n)}$, $D(g)=k$, and $D(hg^{-1})=n$.
Therefore $g$, $h$ satisfy the conditions
of Theorem~\ref{main}(b).~\qed \medskip

     The proof of Theorem~\ref{main} proceeds by cases,
depending mainly on the relative sizes of $k$ and $n$.  We
start with the cases where $n$ is small.  We first require a
lemma.

\begin{lemma} \label{powers}
$D(\f^p)=pk+k_0$.
\end{lemma}

\proof If $p\ge3$ then the result follows from
\cite[Th.\,6]{cam}, while if $p=2$ and $k$ is even then it
follows from \cite[Lemma~1]{cam}.  If $p=2$ and $k$ is odd
then by an explicit calculation we get
$\f^2(x)=x+(r_1r_2+r_1^3)x^4+O(x^5)$ if $k=1$, and
$\f^2(x)=x+r_kr_{k+1}x^{2k+2}+O(x^{2k+3})$ if $k\ge3$,
which implies the result. \qed

\begin{case} \label{case1}
Theorem~\ref{main} holds if $k\le n\le k+k_0$.
\end{case}

\proof For $n$ in this range we have $e(k,n)=k+k_0-n$.
By Lemma~\ref{powers} we have
\begin{equation}
D(\f^p)=D(\g^p)=pk+k_0=n+(p-1)k+e(k,n).
\end{equation}
It follows that $D(\g^p\f^{-p})\ge n+(p-1)k+e(k,n)$, which
proves Theorem~\ref{main}(a).
To prove Theorem~\ref{main}(b) we first consider the
case $n=k+k_0$.  Set
$f(x)=x+x^{k+1}$, $u(x)=x+x^{k+k_0+1}$, and $g=uf$, so that
$D(gf^{-1})=D(u)=k+k_0$.  If $p\mid k$ then $k_0=0$ and $u=f$,
so by \cite[Lemma~1]{cam} we have $D(g^pf^{-p})=D(f^p)=pk$.
If $p\ge3$ and $p\nmid k$ then by the third paragraph in the
proof of \cite[Th.\,6]{cam} we have $D(g^pf^{-p})=pk+k_0$.
If $p=2$ and $n$ is odd then by the explicit computations
in the proof of Lemma~\ref{powers} we get
$D(g^2f^{-2})=2k+1=pk+k_0$.  Thus Theorem~\ref{main}(b) holds
when $n=k+k_0$.  It follows from Lemma~\ref{shift} that
Theorem~\ref{main}(b) also holds for $k\le n<k+k_0$. \qed
\medskip

     We next consider the cases where
$n\ge(p-1)k+p$.  We will need the following
basic result, which is proved in \cite[Prop.\,1]{cam}.

\begin{lemma} \label{basic}
Let $R$ be an integral domain of characteristic $p$,
let $f,g\in\N(R)$, and let $[f,g]=f^{-1}g^{-1}fg$ denote the
commutator of $f$ with $g$.  Then $D([f,g])\ge D(f)+D(g)$, with
equality if and only if $D(f)\not\equiv D(g)\pmod{p}$.
\end{lemma}

     Recall that $\u_1=\g\f^{-1}\in\N(\R)$, and define
$\u_2,\dots,\u_{p}$ inductively by setting
$\u_{i+1}=[\u_i,\f]$.  By Lemma~\ref{basic}
we have $D(\u_{i+1})\ge D(\u_i)+D(\f)$; since $D(\u_1)=n$
this implies $D(\u_i)\geq n+(i-1)k$ for $1\le i\le p$.  It
follows that $D([\u_i,\u_j])\geq 2n+k$ for $1\le i,j\le p$.
Let $\N=\N(\R)$, let $\Nbar=\N/\N_{n+(p-1)k+p}$, and let
$\fbar$, $\gbar$, $\ubar_i$ denote the images of $\f$, $\g$,
$\u_i$ in $\Nbar$.  Since $n\ge(p-1)k+p$ we have
$2n+k\ge n+(p-1)k+p$.  Therefore
$\ubar_i$ commutes with $\ubar_j$ for
all $1\leq i,j\leq p$.  Using the formula
$\u_i\f=\f\u_i\u_{i+1}$ we get
\begin{equation}
\gbar^p=\fbar^p\ubar_1^{C(p,1)}\ubar_2^{C(p,2)}\dots \ubar_p^{C(p,p)},
\end{equation}
where $C(p,i)=p!/i!(p-i)!$ is the binomial coefficient.
It follows from Lemma~\ref{powers} that
\begin{equation}
D(\u_i^p)\ge pn\ge n+(p-1)k+p.
\end{equation}
Since $p\mid C(p,i)$ for $1\leq i\leq p-1$, this implies
$\ubar_i^{C(p,i)}=1$, and hence $\gbar^p=\fbar^p\ubar_p$.
To prove Theorem~\ref{main}
for $n\ge(p-1)k+p$ it suffices to show that
$D(\u_p)\ge n+(p-1)k+e(k,n)$, and that there is a
specialization $\sigma:\R\ra\Fp$ such that
$D(\u_p^{\sigma})=n+(p-1)k+e(k,n)$.

\begin{case} \label{case2}
Theorem~\ref{main} holds if $n\ge(p-1)k+p$ and $p\mid k$.
\end{case}

\proof If $p\nmid n$ then an inductive argument based on
Lemma~\ref{basic} shows that $D(\u_p)=n+(p-1)k$.
If $p\mid n$ then it follows from Lemma~\ref{basic}
that $D(\u_2)\ge n+k+1$, and hence that $D(\u_p)\ge n+(p-1)k+1$.
This proves Theorem~\ref{main}(a).  If $p\nmid n$ let
$\sigma$ be a specialization such that
$D(\f^{\sigma})=k$ and $D(\u_1^{\sigma})=n$.  Then by
Lemma~\ref{basic} we have $D(\u_p)=n+(p-1)k$, which proves
Theorem~\ref{main}(b) in this case.  If $p\mid n$ let
$\sigma$ be a specialization such that $D(\f^{\sigma})=k$
and $D(\u_1^{\sigma})=n+1$.  Using Lemma~\ref{basic} we get
\begin{align}
D(\u_p^{\sigma})&=n+(p-1)k+1 \\
&=n+(p-1)k+e(k,n) \\
&=(n+1)+(p-1)k+e(k,n+1).
\end{align}
Theorem~\ref{main}(b) now follows from Lemma~\ref{shift}.
\qed \medskip

     Since $\u_1$ and $\f$ are generic one might expect that
$D(\u_{i+1})=D(\u_i)+D(\f)+1$ when $D(\u_i)\equiv
D(\f)\pmod{p}$.  In fact this is not always the case: There
are instances where $D(\u_{i+1})>D(\u_i)+D(\f)+1$.  To
compute $D(\u_p)$ when $p\nmid k$, we introduce
a doubly-indexed sequence $(c_{ij})$ which is
closely related to the coefficients of $\u_h(x)$.
We retain the variables $r_k,r_{k+1},r_{k+2},\dots$
and introduce a new variable $K$.  For $i,j\ge0$
we define $c_{ij}\in\Z[K,r_k,r_{k+1},\dots]$ using the
difference equation
\begin{equation} \label{cdiff}
c_{ij}=\sum_{t=0}^j\,\left((i-2)K+n+2t-j\right)
r_{k+j-t}c_{i-1,t}
\end{equation}
for $i\ge1$, $j\ge0$, and the initial conditions
\begin{equation} \label{cinit}
c_{0j}=
\begin{cases}
1&\text{if $j=0$}, \\
0&\text{if $j\ge1$}.
\end{cases}
\end{equation}
Let $S=\Q(K)[r_k^{-1},r_k,r_{k+1},r_{k+2},\dots]$, and for
$a\in\Z$, $i\ge0$ define
\begin{equation}
P_a(i)=\prod_{h=1}^i\,((h-2)K+n+a)\in\Z[K].
\end{equation}

\begin{lemma} \label{coeffs}
There are $\phi_{jab}\in S$ such that for all $i,j\ge0$ we have
\begin{equation} \label{csum}
c_{ij}=r_k^i\cdot\sum_{a=0}^{j}\sum_{b=0}^{j}\,
\phi_{jab}P_a(i+b).
\end{equation}
\end{lemma}

\proof We use induction on $j$.  It follows from
(\ref{cdiff}) and (\ref{cinit}) that $c_{i0}=r_k^iP_0(i)$.
Thus by setting $\phi_{000}=1$ we get the lemma in the case
$j=0$.  Let $j\ge1$ and
assume that the lemma holds for all $c_{it}$ with  $i\ge0$
and $0\le t<j$.  Then for $i\ge1$ and $0\le t<j$ we have
\begin{equation} \label{previous}
c_{i-1,t}=r_k^{i-1}\cdot\sum_{a=0}^t\sum_{b=0}^t\,
\phi_{tab}P_a(i-1+b),
\end{equation}
with $\phi_{tab}\in S$.  Since
\begin{multline}
((i-2)K+n+2t-j)P_a(i-1+b)= \\
P_a(i+b)+(-bK+2t-j-a)P_a(i-1+b),
\end{multline}
by substituting (\ref{previous}) into the difference equation
(\ref{cdiff}) we get
\begin{equation} \label{3sum}
c_{ij}=\left((i-2)K+n+j\right)r_kc_{i-1,j}+
r_k^{i-1}\cdot\sum_{t=0}^{j-1}\sum_{a=0}^t\sum_{b=0}^{t+1}\,
\psi_{tab}P_a(i-1+b),
\end{equation}
where
\begin{equation}
\psi_{tab}=
\begin{cases}
(2t-j-a)r_{k+j-t}\phi_{ta0}&\text{if $b=0$}, \\
(-bK+2t-j-a)r_{k+j-t}\phi_{tab}+
r_{k+j-t}\phi_{t,a,b-1}&\text{if $1\le b\le t$}, \\
r_{k+j-t}\phi_{tat}&\text{if $b=t+1$.}
\end{cases}
\end{equation}

     The general solution to (\ref{3sum}) as a difference
equation in $i$ is
\begin{equation} \label{general}
c_{ij}=\alpha r_k^{i}P_j(i)+
r_k^{i-1}\cdot\sum_{t=0}^{j-1}\sum_{a=0}^t\sum_{b=0}^{t+1}\:
\frac{\psi_{tab}}{bK+a-j}P_a(i+b),
\end{equation}
with $\alpha$ arbitrary.  For $0\le a<j$, $0\le b\le j$ set
\begin{equation}
\phi_{jab}=r_k^{-1}\cdot\sum_{t=c}^{j-1}\:\frac{\psi_{tab}}{bK+a-j},
\end{equation}
where $c=\max\{a,b-1\}$.  Then for $0\le a<b\le j$ we have
\begin{equation} \label{phijab1}
\phi_{jab}=\frac{\phi_{b-1,a,b-1}}{bK+a-j}r_{k+j-b+1}r_k^{-1}+
\sum_{t=b}^{j-1}\,\frac{(-bK+2t-j-a)\phi_{tab}+
\phi_{t,a,b-1}}{bK+a-j}r_{k+j-t}r_k^{-1},
\end{equation}
for $1\le b\le a<j$ we have
\begin{equation} \label{phijab2}
\phi_{jab}=\sum_{t=a}^{j-1}\,\frac{(-bK+2t-j-a)\phi_{tab}+
\phi_{t,a,b-1}}{bK+a-j}r_{k+j-t}r_k^{-1},
\end{equation}
and for $0\le a<j$, $b=0$ we have
\begin{equation} \label{phija0}
\phi_{ja0}=\sum_{t=0}^{j-1}\,\frac{(2t-j-a)\phi_{ta0}}{a-j}r_{k+j-t}r_k^{-1}.
\end{equation}
It follows that (\ref{general}) can be rewritten as
\begin{equation} \label{cij}
c_{ij}=\phi_{jj0}r_k^iP_j(i)+r_k^i\cdot\sum_{a=0}^{j-1}\sum_{b=0}^j
\phi_{jab}P_a(i+b),
\end{equation}
where the value of $\phi_{jj0}=\alpha$ is determined
by the initial conditions (\ref{cinit}) to be
\begin{equation} \label{phijj0}
\phi_{jj0}=-\sum_{a=0}^{j-1}\sum_{b=0}^j\,\phi_{jab}P_a(b).
\end{equation}
Finally, set $\phi_{jjb}=0$ for $1\le b\le j$.  Then $c_{ij}$
is given by (\ref{csum}).
Since $\phi_{tab}\in S$ for $0\le t<j$, it follows from
(\ref{phijab1}), (\ref{phijab2}), and (\ref{phija0}) that
$\phi_{jab}\in S$ for $0\le a<j$ and $0\le b\le j$.
Hence by (\ref{phijj0}) we have $\phi_{jj0}\in S$ as well.
Thus all the coefficients $\phi_{jab}$ in (\ref{csum}) lie
in $S$, so the lemma holds for $j$. \qed \medskip

     In the proof of Lemma~\ref{coeffs} we define $\phi_{jab}$
for all $(j,a,b)$ such that $0\le a,b\le j$.
Set $\phi_{jab}=0$ for all other integer values of $j,a,b$.
Then by (\ref{phijab1}), (\ref{phijab2}), and (\ref{phija0})
we have
\begin{equation}
\sum_{t=a}^{j}\,(bK+j+a-2t)r_{k+j-t}\phi_{tab}=\sum_{t=a}^{j-1}\,
r_{k+j-t}\phi_{t,a,b-1}
\end{equation}
for $j\ge a\ge0$, $b\ge0$.  Shifting $j$ and $t$ by $a$ gives
\begin{equation} \label{phidiff}
\sum_{t=0}^{j}\:(bK+j-2t)r_{k+j-t}\phi_{t+a,a,b}=\sum_{t=0}^{j-1}\,
r_{k+j-t}\phi_{t+a,a,b-1}
\end{equation}
for $j,a,b$ nonnegative.

     From (\ref{phijab1}) and (\ref{phijab2}) we would expect
the denominator of $\phi_{jab}$ to contain the factor
$bK+a-j$.  Surprisingly, this factor is not present unless
$b\mid a-j$.  For $l,m\ge0$ define a subring of $S$,
\begin{equation}
S_{lm}=\Z\left[\frac{1}{l!},K,\frac{1}{K-1},\frac{1}{K-2},\dots,
\frac{1}{K-m},r_k^{-1},r_k,r_{k+1},r_{k+2},\dots\right].
\end{equation}
If $l\le l'$ and $m\le m'$ then clearly $S_{lm}\subset S_{l'm'}$.
For each $(j,a,b)$ we will find $(l,m)$ such that
$\phi_{jab}\in S_{lm}$.  To accomplish this we fix $a$ and find
the generating function for $(\phi_{jab})$.

\begin{prop} \label{Sj}
Let $a\ge0$.
The ordinary generating function for $(\phi_{jab})_{j,b\ge0}$ is
\begin{equation} \label{solution}
F_a(x,y)=\sum_{j=0}^{\infty}\sum_{b=0}^{\infty}\,\phi_{jab}x^jy^b
=\phi_{aa0}r_k^{-1}x^a\alpha(x)e^{\omega(x)y},
\end{equation}
where
$\alpha(x)=r_k+r_{k+1}x+r_{k+2}x^2+\cdots$,
$\dst\frac{r_k}{\alpha(x)}=1+q_1x+q_2x^2+\cdots$, and
\begin{align}
\omega(x)&=\dst\sum_{t=1}^{\infty}\,\frac{q_t}{t-K}\,x^t.
\end{align}
\end{prop}

\proof When $a=0$ the difference equation (\ref{phidiff})
is equivalent to the partial differential equation
\begin{equation} \label{PDE1}
\left(Ky\frac{\partial}{\partial y}+x\frac{\partial}{\partial x}
\right)\left(\strut F_0(x,y)\alpha(x)\right)-
2x\alpha(x)\frac{\partial}{\partial x}F_0(x,y)=yF_0(x,y)\left(
\alpha(x)-r_k\right).
\end{equation}
In addition, since
$\phi_{000}=1$ and $\phi_{00b}=0$ for $b\ge1$, we have
the boundary condition $F_0(0,y)=1$.
The solution to this boundary value problem is easily determined
by the method of characteristics to be the function given in
(\ref{solution}), with $a=0$.  Alternatively, one can use the
formula
\begin{equation}
x\omega'(x)-K\omega(x)=\frac{r_k}{\alpha(x)}-1
\end{equation}
to check directly that this function is the unique solution.
For $a\ge1$, the sequence $(\phi_{j+a,a,b})_{j,b\ge0}$
satisfies the same linear difference equation (\ref{phidiff})
as $(\phi_{j0b})_{j,b\ge0}$, and we have $\phi_{aab}=0$
for $b\ge1$.  Therefore the generating function for
$(\phi_{j+a,a,b})_{j,b\ge0}$ is $\phi_{aa0}F_0(x,y)$.  Since
$\phi_{jab}=0$ for $j<a$, this implies
that the generating function for $(\phi_{jab})_{j,b\ge0}$ is
$F_a(x,y)=\phi_{aa0}x^aF_0(x,y)$.~\qed

\begin{cor} \label{Sm}
Let $j,a,b$ be nonnegative.  If $a+b>j$ then $\phi_{jab}=0$.
If $a+b\le j$ then $\phi_{jab}\in S_{lm}$, where
$l=\max\{a,b\}$, $m=a$ if $b=0$, and $m=\max\{a,j+1-a-b\}$
if $b\ge1$.
\end{cor}

\proof We first consider the case $a=0$.  By
Proposition~\ref{Sj} we see that $\phi_{j0b}$ is the
coefficient of $x^j$ in
\begin{equation} \label{coeff}
r_k^{-1}\alpha(x)\cdot\frac{1}{b!}\cdot\omega(x)^b=
\frac{r_k^{-1}x^b}{b!}\,
(r_k+r_{k+1}x+\cdots)\cdot\left(\sum_{t=1}^{\infty}\,
\frac{q_t}{t-K}\,x^{t-1}\right)^b.
\end{equation}
It follows easily from this observation that $\phi_{j00}$
lies in $S_{00}$,
$\phi_{j0b}$ lies in $S_{b,j+1-b}$ for $1\le b\le j$, and
$\phi_{j0b}=0$ for $b>j$.  We next use induction to show that
$\phi_{jj0}\in S_{jj}$ for $j\ge0$.  Since
$\phi_{000}=1$, the case $j=0$ is clear.  Let $j\ge1$ and
suppose $\phi_{aa0}\in S_{aa}$ for all $0\le a<j$.  By
Proposition~\ref{Sj} we have
$\phi_{jab}=\phi_{aa0}\phi_{j-a,0,b}$, so
(\ref{phijj0}) can be rewritten as
\begin{equation} \label{rewrite}
\phi_{jj0}=-\sum_{a=0}^{j-1}\sum_{b=0}^j\,\phi_{aa0}
\phi_{j-a,0,b}P_a(b).
\end{equation}
By the inductive assumption we have $\phi_{aa0}\in S_{aa}$,
and it follows from the first case that
$\phi_{j-a,0,b}\in S_{jj}$ for $0\le a\le j-1$, $0\le b\le j$.
Therefore all the terms of (\ref{rewrite}) are in $S_{jj}$,
so we get $\phi_{jj0}\in S_{jj}$.  The general case of the
corollary now follows from the formula
$\phi_{jab}=\phi_{aa0}\phi_{j-a,0,b}$. \qed

\begin{cor} \label{nonzero}
View $\phi_{jab}$ as a rational function of $K$.  Then
for each $j\ge1$, $\phi_{j01}$ has a simple pole at $K=j$
with residue $-q_j$.
\end{cor}

\proof Using (\ref{coeff}) we get
\begin{equation}
\phi_{j01}=\frac{r_k^{-1}r_{k+j-1}q_1}{1-K}+\frac{r_k^{-1}r_{k+j-2}q_2}{2-K}+\dots
+\frac{r_k^{-1}r_kq_j}{j-K}.
\end{equation}
The corollary now follows from the fact that $q_j\not=0$.
\qed \medskip

     We now use Corollaries \ref{Sm} and \ref{nonzero} to
prove Theorem~\ref{main} in the cases where $p\nmid k$ and
$n\ge(p-1)k+p$.
Set $e=e(k,n)$ and let $A_h\in \M_{(e+1)\times(e+1)}(\R)$ be
the upper triangular matrix whose $(i,j)$ entry for
$0\le i\le j\le e(k,n)$ is
\begin{equation} \label{ahij}
a_{hij}=((h-2)k+n+2i-j)r_{k+j-i}.
\end{equation}
Then we have
\begin{equation} \label{matrix}
A_h=\left[\begin{array}{ccccc}
n_hr_k&(n_h-1)r_{k+1}&(n_h-2)r_{k+2}&\dots&(n_h-e)r_{k+e}\\
0&(n_h+1)r_k&n_hr_{k+1}&\dots&(n_h-e+2)r_{k+e-1}\\
0&0&(n_h+2)r_k&\dots&(n_h-e+4)r_{k+e-2}\\
\vdots&\vdots&\vdots&&\vdots \\
0&0&0&\dots&(n_h+e)r_k
\end{array}\right],
\end{equation}
where $n_h=(h-2)k+n$.  For $1\le h\le p$ let
$\v_h\in\R^{e+1}$ be the row vector whose entries are the
coefficients of $x^{(h-1)k+n+1},x^{(h-1)k+n+2},\dots,
x^{(h-1)k+n+e+1}$ in $\u_h(x)$.

\begin{lemma} \label{linear}
If $e(k,n)<k$ then $\v_{h+1}=\v_{h}A_h$.
\end{lemma}

\proof Write $\f(x)=x+x^{k+1}\alpha(x)$ and
$\u_h(x)=x+x^{(h-1)k+n+1}\beta(x)$.  Then we have the following
expansions modulo $x^{(h+1)k+n+1}$:
\begin{multline} \label{uhg}
\u_h(\f(x))\equiv x+x^{k+1}\alpha(x)+x^{(h-1)k+n+1}\beta(x)+ \\
((h-1)k+n+1)x^{hk+n+1}\alpha(x)\beta(x)+
x^{hk+n+2}\alpha(x)\beta'(x)
\end{multline}
\begin{multline} \label{guh}
\f(\u_h(x))\equiv x+x^{k+1}\alpha(x)+x^{(h-1)k+n+1}\beta(x)+ \\
(k+1)x^{hk+n+1}\alpha(x)\beta(x)+x^{hk+n+2}\alpha'(x)\beta(x)
\end{multline}
\begin{multline} \label{com}
\![\u_h,\f](x)\equiv x+((h-2)k+n)x^{hk+n+1}\alpha(x)\beta(x)+ \\
x^{hk+n+2}(\alpha(x)\beta'(x)-\alpha'(x)\beta(x)).
\end{multline}
We have $\alpha(x)=r_k+r_{k+1}x+r_{k+2}x^2+\cdots$, and we
can write $\beta(x)=t_0+t_1x+t_2x^2+\cdots$ with
$t_i\in\R$.  Since $e(k,n)\le k-1$, it follows from
(\ref{com}) that for
$0\le j\le e(k,n)$ the coefficient of $x^{hk+n+j+1}$ in
$\u_{h+1}(x)=[\u_h,\f](x)$ is
\begin{equation}
\sum_{i=0}^j\,((h-2)k+n+2i-j)r_{k+j-i}t_i.
\end{equation}
Comparing this expression with (\ref{matrix}) gives the
lemma. \qed

\begin{case} \label{case3}
Theorem~\ref{main} holds if $n\ge(p-1)k+p$, $p\nmid k$,
and $e(k,n)<k$.
\end{case}

\proof For $h\ge1$ define a
matrix $\Pi_h\in\M_{(e+1)\times(e+1)}(\R)$ by setting
$\Pi_h=A_1A_2\dots A_h$.  It follows from Lemma~\ref{linear}
that $\v_p=\v_1\Pi_{p-1}$, where
$\v_1=(s_n,s_{n+1},\dots,s_{n+e})$ has entries which are
independent variables in $\R$ which don't occur in $\Pi_{p-1}$.
Thus to prove
Theorem~\ref{main} in this case it suffices to show that
the first $e(k,n)$ columns of $\Pi_{p-1}$ are all zero, and that
there is a specialization $\sigma:\R\ra\Fp$ which maps
the last column of $\Pi_{p-1}$ to a nonzero element of
$\Fp^{e+1}$.  We indicate the dependence of
$A_h$ on $n$ by writing $A_h=A_h(n)$ and
${a_{hij}=a_{hij}(n)}$.  We also let $\pi_{hij}=\pi_{hij}(n)$
denote the $(i,j)$ entry of $\Pi_h$.  If $0\le i\le j\le e(k,n)$
then $e(k,n+i)\ge e(k,n)-i$, so
$a_{h,0,j-i}(n+i)$ is defined.  By (\ref{ahij}) we have
$a_{hij}(n)=a_{h,0,j-i}(n+i)$ for $0\le i\le j\le e$,
and an inductive argument shows then that
${\pi_{hij}(n)}={\pi_{h,0,j-i}(n+i)}$.  If $j<e(k,n)$ then
$j-i<e(k,n+i)$.  Therefore it will suffice to prove the
following statements for all $n\ge(p-1)k+p$:
\begin{align}
\pi_{p-1,0,j}(n)&=0\text{ for all $0\le j<e(k,n)$,}
\label{zero} \\
\sigma(\pi_{p-1,i,e}(n))&\not=0\text{ for some $0\le i\le e(k,n)$ and
some specialization $\sigma$.}
\label{lastcol}
\end{align}

     There is a natural map $\rho:\Z[K,r_k,r_{k+1},\ldots]\ra\R$
which takes $K$ to the image of $k$ in $\R$.  We denote this
map by $x\mapsto \overline{x}$.  Since $\Pi_h=\Pi_{h-1}A_h$,
the sequence $(\pi_{h0j})$ satisfies the difference equation
\begin{equation} \label{pidiff}
\pi_{h0j}=\sum_{t=0}^j\,\left((h-2)k+n+2t-j\right)
r_{k+j-t}\pi_{h-1,0,t}
\end{equation}
for $h\ge1$, and the initial conditions
\begin{equation} \label{piinit}
\pi_{00j}=
\begin{cases}
1&\text{if $j=0$}, \\
0&\text{if $j\ge1$}.
\end{cases}
\end{equation}
Comparing (\ref{pidiff}) and (\ref{piinit}) with (\ref{cdiff})
and (\ref{cinit}), we see that $\pi_{h0j}=\cbar_{hj}$
for all $h,j$ such that $h\ge0$ and $0\le j\le e(k,n)$.  For
$a\in\Z$, $h\ge0$ let $\Pbar_a(h)=\prod_{i=1}^h\,((i-2)k+n+a)$
denote the image of $P_a(h)$ in $\R$.  Since $p\nmid k$ we have
$\Pbar_a(p-1+b)=0$ for $b\ge1$, and $\Pbar_a(p-1)\not=0$ if and
only if $n\equiv2k-a\pmod{p}$.  Suppose $0\le j<e(k,n)$.  Then
by Corollary~\ref{Sm} we have $\phi_{jab}\in S_{p-1,k_0-1}$
for all $a,b\ge0$.  Since $k_0<p$ the reduction map $\rho$
extends to a map $\tilde{\rho}:S_{p-1,k_0-1}\ra\R[r_k^{-1}]$.
Applying $\tilde{\rho}$ to (\ref{csum}) we get
\begin{equation} \label{zeroes}
\pi_{p-1,0,j}=r_k^{p-1}\cdot\sum_{a=0}^j\sum_{b=0}^j\,
\phibar_{jab}\Pbar_a(p-1+b).
\end{equation}
The terms $\Pbar_a(p-1+b)$ in the sum are all zero, except
those with $b=0$ and $n\equiv2k-a\pmod{p}$.  In this case we
would have $a=e(k,n)$, which contradicts the assumption
$j<e(k,n)$.  Thus $\pi_{p-1,0,j}=0$ for $0\le j<e(k,n)$, which
proves (\ref{zero}).

     To prove (\ref{lastcol}), we first observe that if
$n\equiv 2k-i\pmod{p}$ for some $0\le i\le k_0$ then $e(k,n)=i$
and $\pi_{p-1,e,e}=\Pbar_e(p-1)r_k^{p-1}\not=0$.  Hence
$\sigma(\pi_{p-1,e,e})=\Pbar_e(p-1)\not=0$ for
any $\sigma:\R\ra\Fp$ such that $\sigma(r_k)=1$.
If $n\not\equiv 2k-i\pmod{p}$ for all $0\le i\le k_0$ then
$e(k,n)=k_0$.  By Lemma~\ref{coeffs} we have
\begin{equation} \label{0esum}
c_{p-1,k_0}=r_k^{p-1}\cdot\sum_{a=0}^{k_0}\sum_{b=0}^{k_0}\,
\phi_{k_0ab}P_a(p-1+b).
\end{equation}
We will show that all but two of the terms in (\ref{0esum})
have image zero under $\tilde{\rho}$.  Let $0\le a,b\le k_0$
be such that $(a,b)\not=(0,1)$ and $(a,b)\not=(k_0,0)$.  Then
$\Pbar_a(p-1+b)=0$, and by Corollary~\ref{Sm} we have
$\phi_{k_0ab}\in S_{p-1,k_0-1}$.
Therefore $\phibar_{k_0ab}=\tilde{\rho}(\phi_{k_0ab})$
is defined and $\phibar_{k_0ab}\Pbar_a(p-1+b)=0$.

     It remains to consider the terms $\phi_{k_001}P_0(p)$
and $\phi_{k_0k_00}P_{k_0}(p-1)$ in (\ref{0esum}).  It follows
from (\ref{phijj0}) and the previous paragraph that
$\phi_{k_0k_00}=-\phi_{k_001}P_0(1)+\gamma$ for some
$\gamma\in S_{p-1,k_0-1}$.  Therefore we have
\begin{equation} \label{twoterms}
\phi_{k_001}P_0(p)+\phi_{k_0k_00}P_{k_0}(p-1)=
\phi_{k_001}(P_0(p)-P_0(1)P_{k_0}(p-1))+\gamma P_{k_0}(p-1).
\end{equation}
We wish to expand (\ref{twoterms}) in powers of $K-k_0$.  Let
$Q=\prod_{i=1}^p\,((i-2)k_0+n)\in\Z$.
Then we have the following expansions modulo $(K-k_0)^2$:
\begin{align} \label{P0}
P_0(p)&\equiv Q+Q\cdot\left(\frac{-1}{-k_0+n}+
\sum_{h=0}^{p-2}\,\frac{h}{hk_0+n}\right)(K-k_0) \\
P_0(1)P_{k_0}(p-1)&\equiv Q+Q\cdot\left(\frac{-1}{-k_0+n}+
\sum_{h=0}^{p-2}\,\frac{h-1}{hk_0+n}\right)(K-k_0) \\
P_0(p)-P_0(1)P_{k_0}(p-1)&\equiv \left(
\sum_{h=0}^{p-2}\,\frac{Q}{hk_0+n}\right)(K-k_0). \label{diff}
\end{align}
Since $p\nmid k_0$ and $n\not\equiv k_0\pmod{p}$, there is a
unique $0\le h_0\le{p-2}$ such that $p\mid h_0k_0+n$.  Then
$Q'=Q/(h_0k_0+n)\in\Z$ is the unique term of the sum in
(\ref{diff})
which is not divisible by $p$.  Since $\Pbar_{k_0}(p-1)=0$,
it follows from Corollary~\ref{nonzero} that the image
of (\ref{twoterms}) in $\R$ is
$-\Qbar'\qbar_{k_0}$.  Therefore by
(\ref{0esum}) we get $\pi_{p-1,0,k_0}=\cbar_{p-1,k_0}=
-\Qbar'r_k^{p-1}\qbar_{k_0}$, with
$\Qbar'\in\Fp^{\times}$.  Let $\sigma:\R\ra\Fp$ be a
specialization such that $\sigma(r_i)=0$ for $k<i<k+k_0$.
Then $\sigma(r_k^{p-1}\qbar_{k_0})=
-\sigma(r_k^{p-2}r_{k+k_0})$, so by choosing $\sigma$ so
that $\sigma(r_k)=\sigma(r_{k+k_0})=1$ we get
$\sigma(\pi_{p-1,0,k_0})=\Qbar'\not=0$.  This proves
(\ref{lastcol}). \qed

\begin{case} \label{case4}
Theorem~\ref{main} holds if $n\ge(p-1)k+p$, $p\nmid k$, and
$e(k,n)\ge k$.
\end{case}

\proof If $e(k,n)\ge k$ then $e(k,n)=k=k_0$.  To compute the
necessary coefficients of $\u_h(x)$ we need to
consider the
expansions (\ref{uhg})--(\ref{com}) modulo $x^{(h+1)k+n+2}$.
In this higher-order expansion (\ref{uhg}) and (\ref{com})
acquire the additional term
\begin{equation}
\binom{(h-1)k+n+1}{2}x^{(h+1)k+n+1}\alpha(x)^2\beta(x).
\end{equation}
To account for this extra term the matrix $A_h$ must be
replaced by the upper triangular matrix
$A_h'\in\M_{(k+1)\times(k+1)}(\R)$ whose
$(i,j)$ entry for $0\le i,j\le k$ is
\begin{equation} \label{h0e}
a_{hij}'=
\begin{cases}
a_{h0k}+\binom{(h-1)k+n+1}{2}r_k^2&\text{if $(i,j)=(0,k)$}, \\
a_{hij}&\text{otherwise.}
\end{cases}
\end{equation}
Consequently, $\Pi_h$ is replaced by $\Pi_h'=A_1'A_2'\dots A_h'$.
An easy computation shows that there is $m\in\Fp$ such that
\begin{equation}
\pi_{p-1,i,j}'=
\begin{cases}
\pi_{p-1,0,k}+mr_k^p&\text{if $(i,j)=(0,k)$}, \\
\pi_{p-1,i,j}&\text{otherwise.}
\end{cases}
\end{equation}
Thus the first $k$ columns of $\Pi_{p-1}'$ are the same as
those of $\Pi_{p-1}$.  It follows by the reasoning in
Case~\ref{case3} that the first $k$ columns
of $\Pi_{p-1}'$ are zero, and hence that Theorem~\ref{main}(a)
holds in this case.  If $n\equiv 2k-i\pmod{p}$ for
some $0\le i\le k$ then $\pi_{p-1,k,k}'=\pi_{p-1,k,k}$, so
there is a specialization $\sigma$ such that
$\sigma(\pi_{p-1,k,k}')=\sigma(\pi_{p-1,k,k})\not=0$ as in
Case~\ref{case3}.
If $n\not\equiv 2k-i\pmod{p}$ for all $0\le i\le k$ then
$\pi_{p-1,0,k}'=-\Qbar'r_k^{p-1}\qbar_k+mr_k^p$.  Let
$\sigma$ be a specialization such that
$\sigma(r_i)=0$ for $k<i<2k$.  As in Case~\ref{case3} we have
$\sigma(r_k^{p-1}\qbar_{k_0})=-\sigma(r_k^{p-2}r_{2k})$,
and hence $\sigma(\pi_{p-1,0,k}')=
\Qbar'\sigma(r_k^{p-2}r_{2k})+m\sigma(r_k^p)$.  By
choosing $\sigma$ so that $\sigma(r_k)=1$ and
$\sigma(r_{2k})$ is either 0 or 1
we get $\sigma(\pi_{p-1,k,k}')\not=0$.
Therefore Theorem~\ref{main}(b) holds in this case. \qed \medskip

     In Case \ref{case1} we proved Theorem~\ref{main} for
$n\le k+k_0$, and in Cases \ref{case2}, \ref{case3},
\ref{case4} we proved Theorem~\ref{main} for $n\ge(p-1)k+p$.
It remains to prove Theorem~\ref{main} for intermediate values
of $n$.  Recall that $\f(x)=x+r_kx^{k+1}+r_{k+1}x^{k+2}+\cdots
\in\N(\R)$.  The map
$T_{\f}:\N(\R)\ra\N(\R)$ defined by $T_{\f}(h(x))=h(\f(x))$
induces a linear transformation on coefficient vectors
$(1,a_1,a_2,\dots)$ of elements of $\N(\R)$ (see
\cite{york}).  This linear transformation can be
represented by right multiplication by an infinite
matrix of the form $I+M$, where $I$ is the
identity and $M$ is an upper triangular matrix whose diagonal
entries are all 0.  The rows and columns of $I$ and $M$ are
indexed by positive integers, and for $1\le i<j$ the
$(i,j)$ entry of $M$ is
\begin{equation} \label{mij}
m_{ij}=\sum_{l_1+\cdots+l_i=j-i}r_{l_1}r_{l_2}\dots r_{l_i},
\end{equation}
where by convention we set $r_0=1$ and $r_l=0$ for $1\le l<k$.
In particular, we have $m_{ij}=0$ for $j-i<k$.
The formula (\ref{mij}) may be rewritten as
\begin{equation} \label{sum}
m_{ij}=\sum\,\binom{i}{n_0,n_k,n_{k+1}\dots,n_{j-i}}
r_k^{n_k}r_{k+1}^{n_{k+1}}\dots r_{j-i}^{n_{j-i}},
\end{equation}
where the sum is taken over nonnegative integers
$n_0,n_k,n_{k+1},\dots,n_{j-i}$ such that
\begin{align}
n_0+n_k+n_{k+1}+\dots+n_{j-i}&=i \\
kn_k+(k+1)n_{k+1}+\dots+(j-i)n_{j-i}&=j-i.
\end{align}
For $h\ge2$ the $(i,j)$ entry of $M^h$ can be expressed in
terms of the $m_{ij}$:
\begin{equation} \label{mijh}
m_{ij}^{(h)}=\sum_{i<b_1<\dots<b_{h-1}<j}m_{ib_1}m_{b_1b_2}
\dots m_{b_{h-1}j}.
\end{equation}
This formula allows us to compute the entries of the
matrix $(I+M)^p=I+M^p$ which represents $\f^p$.
We are particularly interested in $m_{1j}^{(p)}$,
which for $j\ge2$ is the coefficient of $x^j$ in $\f^p(x)$.

\begin{lemma} \label{modp}
Let $d,i$ be positive integers such that $d<pk$ and
$i\ge(d+1-k)/k$.  Then $m_{i,i+d}=m_{i+p,i+p+d}$.
\end{lemma}

\proof Note that if $n_k,n_{k+1},\dots,n_d$ are nonnegative
integers such that
\begin{equation} \label{knkd}
kn_k+(k+1)n_{k+1}+\dots+dn_d=d
\end{equation}
then we have $n_k+n_{k+1}+\dots+n_d\le d/k$.  Since
$i\ge(d+1-k)/k$ this implies
\begin{equation}
i-(n_k+n_{k+1}+\dots+n_d)\ge\frac{d+1-k}{k}-\frac{d}{k}>-1.
\end{equation}
It follows that $n_0=i-(n_k+n_{k+1}+\dots+n_d)$ is
nonnegative.  Hence by (\ref{sum}) we have 
\begin{equation} \label{miid}
m_{i,i+d}=\sum\,\binom{i}{n_0,n_k,n_{k+1},\dots,n_d}
r_k^{n_k}r_{k+1}^{n_{k+1}}\dots r_d^{n_d},
\end{equation}
where the sum is taken over nonnegative
$n_k,n_{k+1},\dots,n_d$ satisfying (\ref{knkd}).
For $k\le j\le d$
we have $jn_j\le d<pk\le pj$, and hence $n_j<p$.  Therefore
in characteristic $p$ the multinomial coefficient
\begin{align}
\binom{i}{n_0,n_k,n_{k+1},\dots,n_d}&=
\frac{i!}{n_0!n_k!n_{k+1}!\dots n_d!} \\[.2cm]
&=\frac{i(i-1)\dots(i-(n_k+n_{k+1}+\dots+n_d-1))}{n_k!n_{k+1}!\dots n_d!}
\end{align}
is unchanged if we replace $i$ by $i+p$ and $n_0$ by $n_0+p$.
It follows that (\ref{miid}) is also unchanged if we replace
$i$ by $i+p$. \qed

\begin{lemma} \label{expand}
Let $0\le t<n$ and $i\ge2$.  Then there is a
polynomial $A$ with coefficients in $\Fp$ such that
\begin{equation}
m_{i,i+n+t}=A(r_k,r_{k+1},\dots,r_{n-1})+ir_{n+t}+
i\cdot\sum_{w=0}^{t-k}\,r_{n+w}m_{i-1,i-1+t-w}.
\end{equation}
\end{lemma}

\proof We can write
\begin{align}
m_{i,i+n+t}&=\sum_{l_1+\cdots+l_i=n+t}r_{l_1}r_{l_2}\dots r_{l_i}
\\
&=A(r_k,r_{k+1},\dots,r_{n-1})+B(r_k,r_{k+1},\dots,r_{n+t}),
\end{align}
where $A$ is the sum of the terms which depend only on
$r_k,r_{k+1},\dots,r_{n-1}$, and $B$ is the sum of the
remaining terms.  Let $r_{l_1}r_{l_2}\dots r_{l_i}$ be a term
of $B$.  Then $l_j\ge n$ for some $j$, so we have
$l_j=n+w$ with $0\le w\le t$.  Furthermore, since $t<n$, we
have $l_h<n$ for all $h\not=j$.  Since there are $i$ possible
values for $j$ we get
\begin{equation}
B(r_k,r_{k+1},\dots,r_{n+t})=ir_{n+t}+i\cdot\sum_{w=0}^{t-1}\,
r_{n+w}m_{i-1,i-1+t-w}.
\end{equation}
Since $m_{i-1,i-1+t-w}=0$ for $w>t-k$, the lemma follows. \qed

\begin{prop} \label{addp}
Let $s\ge0$ satisfy $n>k+s$ and $pk>k+s$.  Then the
coefficient of $x^{1+s+n+(p-1)k}$ in $\f^p(x)$ can be written
uniquely in the form
\begin{equation} \label{Ens}
m_{1,1+s+n+(p-1)k}^{(p)}=C_{ns}+E_{ns}^{(0)}r_n+E_{ns}^{(1)}r_{n+1}+
\dots+E_{ns}^{(s)}r_{n+s},
\end{equation}
with $C_{ns}\in\Fp[r_k,r_{k+1},\dots,r_{n-1}]$ and
$E_{ns}^{(w)}\in\Fp[r_k,r_{k+1},\dots,r_{k+s-w}]$ for
$0\leq w\leq s$.  Furthermore, we have
$E_{n+p,s}^{(w)}=E_{ns}^{(w)}$.
\end{prop}

\proof Set $b_0=1$ and $b_p=1+s+n+(p-1)k$.  By (\ref{mijh})
we have 
\begin{equation} \label{mp}
m_{1,1+s+n+(p-1)k}^{(p)}=
\sum_{b_0<b_1<\dots<b_p}m_{b_0b_1}m_{b_1b_2}\dots m_{b_{p-1}b_p}.
\end{equation}
To prove the first statement it suffices to show that each
term in the sum (\ref{mp}) can be expressed in the form of
(\ref{Ens}), i.\,e.,
\begin{equation} \label{term}
m_{b_0b_1}m_{b_1b_2}\dots m_{b_{p-1}b_p}=
c+e^{(0)}r_n+e^{(1)}r_{n+1}+\dots+e^{(s)}r_{n+s},
\end{equation}
with $c\in\Fp[r_k,r_{k+1},\dots,r_{n-1}]$ and
$e^{(w)}\in\Fp[r_k,r_{k+1},\dots,r_{k+s-w}]$.  If
$m_{b_0b_1}m_{b_1b_2}\dots m_{b_{p-1}b_p}$ lies in
$\Fp[r_k,r_{k+1},\dots,r_{n-1}]$ then this is clear.
If $m_{b_0b_1}m_{b_1b_2}\dots m_{b_{p-1}b_p}$ depends on
$r_h$ for some $h\ge n$ then $b_i-b_{i-1}\ge k$
for all $1\le i\le p$ and $b_{j}-b_{j-1}=n+t$ for some
$1\le j\le p$ and $t\ge0$.  For $i\not=j$
we have $b_{i}-b_{i-1}\le k+s-t$.  Therefore $t\le s$ and
$m_{b_{i-1}b_i}\in\Fp[r_k,r_{k+1},\dots,r_{k+s-t}]$. 
If $j=1$ then $b_{j-1}=b_0=1$, and hence
$m_{b_0b_1}=m_{1,1+n+t}=r_{n+t}$.  In this case
$m_{b_0b_1}m_{b_1b_2}\dots m_{b_{p-1}b_p}$ can be written in
the form (\ref{term}) with $c=0$, $e^{(t)}=m_{b_1b_2}\dots
m_{b_{p-1}b_p}$, and $e^{(w)}=0$ for all $w\not=t$.
If $j\ge2$ then $b_{j-1}\ge2$, and hence $m_{b_{j-1}b_j}$
is given by Lemma~\ref{expand}.  Therefore
$m_{b_0b_1}m_{b_1b_2}\dots m_{b_{p-1}b_p}$ can be written in
the form (\ref{term}), with $e^{(w)}=0$ for all $w\not=t$
with $w>t-k$.  It follows that $m_{1,1+s+n+(p-1)k}^{(p)}$ can
be written in the form (\ref{Ens}).  The fact that $C_{ns}$
and $E_{ns}^{(w)}$ are uniquely determined follows from the
assumption $n>k+s$.

     To prove the last statement we observe
that there is a one-to-one correspondence between the
terms of (\ref{mp}) which don't lie in
$\Fp[r_k,r_{k+1},\dots,r_{n-1}]$ and the terms of the
corresponding expansion of $m_{1,1+s+n+p+(p-1)k}^{(p)}$ which don't
lie in $\Fp[r_k,r_{k+1},\dots,r_{n+p-1}]$.  This
correspondence is given by
\begin{equation}
m_{b_0b_1}\dots m_{b_{j-1}b_j}\dots m_{b_{p-1}b_p}
\longleftrightarrow m_{b_0b_1}\dots m_{b_{j-1},b_j+p}\dots
m_{b_{p-1}+p,b_p+p},
\end{equation}
where $j$ is determined by the condition
$b_j-b_{j-1}\ge n$.  Set $b_j-b_{j-1}=n+t$ as above.  It
follows from
Lemma~\ref{expand} that for $0\le w\le t-k$ the coefficient
of $r_{n+w}$ in $m_{b_{j-1}b_j}$ and the coefficient of
$r_{n+p+w}$ in $m_{b_{j-1},b_j+p}$ are both equal to
$b_{j-1}m_{b_{j-1}-1,b_{j-1}-1+t-w}$, and that the
coefficient of $r_{n+t}$ in $m_{b_{j-1}b_j}$ and the
coefficient of $r_{n+p+t}$ in $m_{b_{j-1},b_j+p}$ are both
equal to $b_{j-1}$.  In addition, since $pk>k+s$, Lemma~\ref{modp}
implies that $m_{b_{i-1}b_i}=m_{b_{i-1}+p,b_i+p}$ for
$i>j$.  It follows that $E_{n+p,s}^{(w)}=E_{ns}^{(w)}$ for
$0\le w\le s$.~\qed

\begin{case} \label{case5}
Theorem~\ref{main} holds if $k+k_0<n<(p-1)k+p$.
\end{case}

\proof Choose $n'$ such that $n'\equiv n\pmod{p}$ and
$n'\ge (p-1)k+p$, and set $e=e(k,n)=e(k,n')$.  It follows
from Theorem~\ref{main}(a) in Cases \ref{case2}, \ref{case3},
\ref{case4} and Remark~\ref{interp} that $E_{n's}^{(w)}=0$
for all $s,w$ such
that $0\le w\le s<e$.  It follows from Proposition~\ref{addp}
that $E_{ns}^{(w)}=0$ for $0\le w\le s<e$.  Therefore by
Remark~\ref{interp} we see that Theorem~\ref{main}(a) holds
for $n$.  It follows from Theorem~\ref{main}(b)
in Cases \ref{case2}, \ref{case3}, \ref{case4} that there
is a specialization $\sigma$ such
that $\sigma(E_{n'e}^{(w)})\not=0$ for some $0\le w\le e$.
Using Proposition~\ref{addp} we get
$\sigma(E_{ne}^{(w)})=\sigma(E_{n'e}^{(w)})\not=0$.
Let $\tau:\R\ra\Fp$ be a specialization such that
$\tau(r_i)=\sigma(r_i)$ for $i\not=n+w$, and
$\tau(r_{n+w})\not=\sigma(r_{n+w})$.  Since
$\tau(E_{ne}^{(w)})=\sigma(E_{ne}^{(w)})\not=0$, it follows
from Proposition~\ref{addp} that
\begin{equation}
D((\f^{\tau})^p(\f^{\sigma})^{-p})=n+(p-1)k+e(k,n).
\end{equation}
Therefore Theorem~\ref{main}(b) holds for $n$. \qed \medskip

     By combining Cases \ref{case1} through \ref{case5} we
conclude that Theorem~\ref{main} holds for all
$n\ge k\ge1$.\medskip

\noindent
{\bf Acknowledgment:} I would like to thank
Professor Charles Leedham-Green for asking whether
Theorem~\ref{main}(b) is true.

\end{document}